\documentclass[dvipsnames, 12pt]{scrartcl}

\usepackage{pgfplots}
\usepackage{subcaption}
\usepackage[english]{babel} 
\usepackage{amsmath}
\usepackage{mathtools}
\usepackage{multirow}
\usepackage{tablefootnote}
\usepackage{bm}
\usepackage[capitalise]{cleveref}
\usepackage[standard]{ntheorem}
\usepackage{tikz}
\usepackage[draft, margin]{fixme}
\fxsetup{theme=color}
\usepackage{graphicx}
\usepackage[algoruled]{algorithm2e}
\usepackage{grffile}
\usepackage[backend=biber,style=numeric,maxbibnames=99]{biblatex}
\usepackage{fancyhdr}
\usepackage{enumitem}
\usepackage{float}
\usepackage{nicefrac}
\usepackage{csquotes}
\usepackage{booktabs}
\usepackage{textcomp}
\usepackage{bbm}
\usepackage[T1]{fontenc}
\usepackage{amssymb}

\FXRegisterAuthor{pg}{apg}{PG}
\FXRegisterAuthor{mf}{amf}{MF}
\FXRegisterAuthor{fk}{afk}{FK}
\FXRegisterAuthor{aa}{aaa}{AA}

\crefname{assumption}{Assumption}{Assumptions}

\crefname{claim}{Claim}{Claims}
\crefname{equation}{}{}  
\Crefname{equation}{Equation}{Equations}

\DeclareMathOperator{\cost}{cost}

\DeclareMathOperator\interior{int}

\newcommand{\Cs}{\mathcal{C}}

\newcommand{\Kappa}{\mathrm{K}}

\newcommand{\lambdamax}{\lambda^{+}(\Acc)}
\newcommand{\lambdamin}{\lambda^{-}(\Acc)}
\newcommand{\AK}{{\Acc,S}}

\newcommand{\1}{{\mathbbm{1}}}

\newcommand{\coresetX}{\tilde{X}}

\newcommand{\coresetOmega}{\tilde{\Omega}}
\newcommand{\coresetC}{\tilde{C}}
\newcommand{\supportX}{X'}

\newcommand{\supportOmega}{\Omega'}

\newcommand{\extension}{g}

\newcommand{\gridX}{X}

\newcommand{\Acc}{\mathcal{A}}
\newcommand{\Ecc}{\mathcal{E}}

\newcommand{\wcaa}{$\mathbb{WCAA}$\xspace}
\newcommand{\wcac}{weight-constrained anisotropic clustering\xspace}

\newcommand{\wcacs}{weight-constrained anisotropic clusterings\xspace}

\newcommand{\DiM}{DiLPM\xspace}

\newcommand{\apd}{\mathcal{P}}
\newcommand{\APD}{\mathbb{APD}}

\newcommand{\apdmatrices}{\mathcal{A}}
\newcommand{\apdmatrix}{A}
\newcommand{\apdweight}{\gamma}
\newcommand{\apdweights}{\Gamma}
\newcommand{\apdsites}{S}
\newcommand{\apds}{s}

\newcommand{\fa}{\mathfrak{A}}
\newcommand{\fx}{\mathfrak{X}}

\DeclarePairedDelimiterX{\norm}[1]{\lVert}{\rVert}{#1}
\DeclarePairedDelimiterX{\set}[1]{\lbrace}{\rbrace}{#1}

\newcommand{\R}{\mathbb{R}}
\newcommand{\N}{\mathbb{N}}
\newcommand{\Z}{\mathbb{Z}}

\makeatletter
\newcommand{\leqnomode}{\tagsleft@true\let\veqno\@@leqno}
\newcommand{\reqnomode}{\tagsleft@false\let\veqno\@@eqno}
\makeatother

\pgfplotsset{compat=newest} 
\bibliography{references}

\begin{document}
	\selectlanguage{english}
	\title{Turning grain maps into diagrams}
	\author{Andreas Alpers, Maximilian Fiedler, \\ Peter Gritzmann, Fabian Klemm}
	
	\publishers{\vspace*{4ex}%
		\normalfont\normalsize%
		\parbox{0.8\linewidth}{%
			\textbf{Abstract.} 
			The present paper studies mathematical models for representing, imaging, and analyzing polycrystalline materials. We introduce various techniques for converting grain maps into diagram or tessellation representations that rely on constrained clustering. In particular, we show how to significantly accelerate the generalized balanced power diagram method from~\cite{Alpers2015} and how to extend it to allow for optimization over all relevant parameters. A comparison of the accuracies of the proposed approaches is given based on a 3D real-world data set of $339\times 339 \times 599$ voxels.
		}	
	}
	
	\maketitle


\section{Introduction}\label{sec:introduction}

A challenging problem in materials science is to describe the structure of polycrystalline materials and understand the processes that govern the evolution of the structure with time, e.g., the grain growth taking place during heating; see, e.g.,~\cite{hansenpoulsen2009,poulsenbook,banhartbook}. Naturally,
geometric concepts play an important role, e.g.,~\cite{model3}, and computations can utilize voxel-based images at different levels of resolution, in terms of 2D~\cite{teferra18,oleg2012} or 3D \emph{grain maps}~\cite{Zhang2018}. However, the data acquisition and analysis tend to be slow, prohibiting the study of real-life dynamics. Moreover, the size of some grains may be comparable to the spatial resolution of the methods, rendering detailed visualization of shape impossible. For these reasons, sparse descriptions involving only a few parameters have been introduced. Such \enquote{compressed segmentations} are also relevant for materials modelling, as the reduced complexity allows for the analysis of much larger (and more representative) samples.

The present paper is concerned with deterministic mathematical models for converting polycrystalline images into diagram representations efficiently.  We perform a \enquote{controlled analysis} based on different types of image discretization,  resolutions, and aspects of data compression. Such discrete images, formally introduced in \cref{sec:diagrams}, will generically be referred to as \enquote{grain scans}. (This term is more general than the standard term grain map from materials science.) The main goal of our analysis is to further substantiate the potential of diagram representations with a special view, however, on later standard applications for solving the inverse problem of grain mapping when only certain measurements such as centroids, volumes, or centralized second-order moments of grains are available. 

Diagrams have been applied for representations of static polycrystalline materials before. The known models, such as \emph{Voronoi} \cite{Xu2009,barker16}, \emph{Voronoi~G} \cite{guilleminot11}, \emph{ Laguerre} \cite{Lyckegaard2011,laguerre2,Telley1992, Altendorf2014,kuehn08}, \emph{ellipsoidal growth} \cite{Teferra2015}, and \emph{generalized balanced power diagrams (GBPD)} \cite{Alpers2015} vary in terms of their accuracy and their complexity. It has been observed in several studies, including \cite{sedivy16, sedivy17,teferra18, Spettl2016, curvedgrainboundaries} that the anisotropic power diagrams produced by the GBPD-algorithm capture the physical principles governing the forming of polycrystals quite well. However, the high accuracy comes at the price of high computational costs, which limits the practical use of the method. In fact, all known approaches are either heuristic in nature or require a significant computational effort. Indeed, there is some need for balancing the computational effort and the desired accuracy for most applications. This problem is aggravated by a strong increase in the size of grain models both experimentally and from a material modelling point of view.

The present paper introduces new and significantly accelerated variants of the GBPD algorithm of \cite{Alpers2015}, which generates GBPDs, a diagram class of large flexibility containing Voronoi and Laguerre tessellations as special cases (see Section~\ref{sec:diagrams}). In particular, we utilize data compression techniques based on \emph{coresets} which allow for a substantial \enquote{thinning} of the data without loosing too much accuracy. In fact, we can compute optimal diagrams based on low-resolution ground set models whose number of points depend only on the number of grains and the desired accuracy. We will also introduce additional techniques resulting in yet sparser image supports, altogether generating a class of new algorithms for computing diagram representations, which we will call \emph{sparse-GBPD} or  \emph{s-GBPD}. As it turns out, suitable parameter settings yield algorithms in s-GBPD, which are dramatically faster than the GBPD algorithm at essentially no loss of accuracy, rendering the new algorithm a viable, practical tool of high accuracy. While s-GBPD is an \enquote{indirect} clustering method based on certain experimentally measurable parameters, we will also introduce a new technique, addressing the best fit problem directly, which allows optimizing over all relevant parameters that govern the weight-balanced anisotropic diagrams. 

The paper is organized as follows. \cref{sec:diagrams}  formalizes the concept of grain scans, introduces diagrams, states the main problems and indicates our main contributions. \cref{sec:proxiess-GBPD} introduces the concept of coresets for speeding up exact methods and gives several examples of how such data compression methods can be realized in our context. In particular, we show how the original GBPD algorithm can be turned into a computationally feasible method for practical 3D problem sizes and reasonable resolutions. \cref{sec:full-model}  introduces a new, direct, and full-parameter best fit model. The results of our empirical study comparing the different methods on available practical grain scan data are stated in \cref{sec:evaluation_polycrystals}. Finally, \cref{sec:conclusion} concludes with remarks that particularly address the issues of how the coreset methods can be extended to situation where instead of full grain scans only certain grain measurements are available.

\section{Turning grain scans into diagrams}\label{sec:diagrams}

Adequate representations of grain scans can be valuable tools for understanding the physical principles of grain growth and the properties of the resulting material. The diagram representations studied in this paper require only a few parameters per grain and are thus much easier to process than voxel-represented grain scans. Further, since diagram representations provide a dissection of the underlying (continuous) space they can even be used to boost the scan resolution. Moreover, they exhibit geometric and combinatorial features of the grains that are not readily available otherwise. Hence, diagram representations may contribute to predicting and even controlling the forming of new materials.

As previous studies \cite{Alpers2015,sedivy16, sedivy17,Spettl2016, curvedgrainboundaries, teferra18} show, only a few characteristic parameters govern decisive properties of grain scans. These parameters may, in practice, be derived from measurements of existing polycrystals or demanded by new materials. Here, we take a more fundamental perspective and assume that grain scans are available at some resolution. A specific focus will lie on the accuracy of fit of the derived diagrams. Also, computations should be fast but not necessarily real-time as we may want to spend additional time improving the accuracy of the representation. 

In this section we will first formalize the concept of a grain scan and introduce the diagrams used throughout the paper. Then we put our studies in the subsequent sections into some broader perspective and state and discuss our main contributions.

\subsubsection*{Grain scans}
 Let $d\in \N$ be the dimension we are working in (e.g., $d=2$ and $d=3$ for the study of planar or 3D samples, respectively), and let $k\in \N$ denote the number of different grains in the given polycrystalline sample. Given a set of points  $X=\{x_1,\dots,x_n\}\subseteq \mathbb{R}^d$ in the sample and labels, i.e., $y_1,\dots,y_n\in[k]:=\{1,\dots,k\},$ we call the labeled data set $\bigl\{(x_j,y_j): j\in [n]\bigr\}$ a \emph{grain scan} or \emph{polycrystal scan}. We will refer to it as $(X,Y)$ with the understanding, that the indices of the labels in the family $Y= \{ y_1,\dots , y_n \}$ correspond to the indices of the points in $X$.
 
Typically our samples (are normalized to) lie in the range $[0,1]^d$, and often (but, as we also work with different sets of compressed data, not exclusively), the points of $X$ are the grid points $(\nicefrac{1}{\rho_1} \cdot \Z \times \dots \times \nicefrac{1}{\rho_d}\cdot \Z)\cap [0,1]^d$ for some \emph{resolution} vector $\rho=(\rho_1,\dots,\rho_d) \in \N^d$. Such special cases can be viewed as digital images, commonly referred to as \emph{grain maps}.  
 
 In this notation, we can identify the (discretization of the) $i$th grain with the set
\begin{equation*}
    G_i= \bigl\{x_j\in X: (x_j,i)\in (X,Y)\bigr\},
\end{equation*}
and, with $\kappa_i:=|G_i|$ we can view the product $\kappa_i \cdot \nicefrac{1}{\rho_1} \cdot \nicefrac{1}{\rho_2} \cdot \ldots \cdot \nicefrac{1}{\rho_d}$ as an approximation of $G_i$'s volume.

\subsubsection*{Diagrams}

Given a grain scan $(X,Y)$ we are interested in a diagram that represents the scan as accurately as possible. In its general form, \emph{diagrams} can be defined by functions $h_1,\ldots,h_k : \R^d \rightarrow \R$. For each $i \in [k]$, the corresponding \emph{cell} $P_i$ of the diagram is given as 
\begin{align*}
	P_i = \bigl\{ x \in \R^d:  h_i(x) \leq   h_l(x) ~~\forall l \in [k]\bigr\}.
\end{align*}
The choice of the functions $h_i$ determines the diagram type; see \cite{Brieden2017} and \cite{GK17} for examples and properties of different classes of diagrams.
The classical \emph{Voronoi diagram} of $k$ distinct points $\apds_1, \ldots, \apds_k$ in $\R^d$, for instance, is obtained by employing the functions $h_i= \norm{x-\apds_i}_{2}$ where $\norm{\quad}_{2}$ denotes the Euclidean norm. In the following we use three types of parameters to controll the properties of the cells:
\begin{itemize}
  \item[$\apdmatrices$:\quad] $k$ positive definite symmetric matrices $\apdmatrix_{1},\dots,\apdmatrix_{k} \in \R^{d\times d}$, collected in the family $\apdmatrices=\{\apdmatrix_{1},\dots,\apdmatrix_{k}\}$. In the special case that each $\apdmatrix_{i}$ is the unit matrix $E_d$ of $\R^{d\times d}$ the family $\apdmatrices$ will be denoted by $\Ecc$; 
  \item[$\apdsites$:\quad] $k$ different \emph{sites} $\apds_1, \ldots, \apds_k \in \R^d$, collected in the set $\apdsites = \{ \apds_1, \ldots, \apds_k\}$;  
 \item[$\Gamma$:\quad] $k$ \emph{sizes} $\apdweight_1,\dots, \apdweight_k \in \R$, collected in the family $\Gamma = \{ \apdweight_1,\dots, \apdweight_k\}$.  
\end{itemize}
The sites specify the general positions of the grains while the matrices in $\apdmatrices$ describe characteristics of their shapes. In fact, each grain is equipped with its own \emph{ellipsoidal norm} $\norm{\quad }_{A_i}$ given through $\apdmatrix_i$ by
\begin{align*}
	\norm{ x }_{A_i} := \sqrt{ x^\top A_i x }, \qquad (x \in \R^d).
\end{align*}
Finally, the sizes $\apdweight_i$ can be used to control the grain volume. 

The \emph{anisotropic power diagram} (APD)
\begin{align*}
\apd= \apd(\apdmatrices,\apdsites,\Gamma) =\{P_1,\ldots,P_k\}
\end{align*}
is then the diagram specified by the functions
\begin{equation*}
    h_i(x) = \norm{x-s_i}_{A_i}^2 + \apdweight_i \qquad (i\in [k]).
\end{equation*}
Let us point out that $\apd$ constitutes a cell decomposition of $\R^d$ and, in particular, two different cells do not share interior points, i.e., $\interior(P_i)\cap \interior(P_\ell)=\emptyset$ for $i\ne \ell$; see \cite{Brieden2017} for a proof (even for more general diagrams).

Clearly, APDs generalize the \emph{Voronoi diagrams} (where $A_i=E_d$ and $\gamma_i=0$ for all $i\in [k]$) and also the \emph{power diagrams} (PD) (where $\apdmatrices=\Ecc$) whose cells are polyhedra. While the Euclidean norm $\norm{\quad }_{E_d}=\norm{\quad }_{2}$ corresponds to the \emph{isotropic} case, the more general ellipsoidal norms allow us to model also \emph{anisotropic} growth. 

\subsubsection*{Representation}

Given a polycrystalline scan, we are interested in an APD which \emph{represents} the scan, i.e., ideally $G_i = X \cap P_i$ for all $i \in [k]$. Of course, in general, we cannot expect a perfect fit. Hence, we resort to some notion of \emph{best fit}.

In order to define the accuracy of an APD $\apd=\{P_1,\ldots,P_k\}$, we consider the associated \emph{classification function} 
\begin{equation*}
 f_\apd:\R^d \rightarrow [k] \cup \{0\}, \qquad x \mapsto 
 \begin{cases}
 i & \text{if $x\in \interior(P_i)$;}\\
 0& \text{otherwise}.
 \end{cases}
\end{equation*}
The function $f_\apd$ maps a point $x \in X$ either to the index of its unique containing cell or (and this is a conservative setting) to $0$ if $x$ lies on the boundary of more than one cell. Using such functions, we can formulate the problem of best-representing polycrystal scans as finding a diagram $\apd$ such that $f_\apd$ minimizes the number of misclassified data points. More precisely, we want to solve the following problem.

\begin{quote}
    $\APD$:\\
    Given $n \in \N$ and an $n$-element labeled data set $(X,Y)$, compute, for $i \in [k]$, matrices $\apdmatrix_i$, sites $\apds_i$ and sizes $\apdweight_i$ such that the \emph{classification error}
\begin{equation*}
		\left|\bigl\{j \in [n]: f_\apd(x_j) \neq y_j\bigr\}\right|
\end{equation*}
of the corresponding diagram $\apd$ is minimized over all feasible choices of ($\apdmatrices,\apdsites,\apdweights$). Moreover, return this classification error.
\end{quote}
For $\APD$, all three characteristics $\apdmatrices$, $\apdsites$, and $\apdweights$ are optimized and constitute the output of the problem, together with the classification error. As this reduces the computational effort, it is often useful to specify some of these characteristics beforehand and optimize the others. Fixing characteristics will be denoted by adding them in brackets. For instance, $\APD(\apdmatrices,\apdsites$) refers to the situation considered in \cite{Alpers2015} that both, $\apdmatrices$ and $\apdsites$, are given and the optimization involves only $\apdweights$. Even more extreme, $\APD(\apdmatrices,\apdsites,\apdweights)$ indicates that all characteristics are given, therefore, the anisotropic power diagram $\apd=\apd(\apdmatrices,\apdsites,\apdweights)$ is completely determined, and only the classification error evoked by $\apd$ is to be computed. 

Since both $d$ and $k$ are fixed, it follows from arguments given in \cite{BBG17}, see also \cite{MaxDiss}, that $\APD$ can be solved by enumerating all relevant APDs in polynomial time. Their number grows, however, exponentially in $k$ and $d$. Hence, such an enumerative approach is computationally infeasible for virtually all practically relevant instances. The known APX-hardness results of \cite{Ben-David2003} and \cite{MaxDiss} for the case that $k$ is part of the input indicate that there is not much hope to find dramatically better algorithms for solving $\APD$ exactly. Accordingly, various heuristics suggest to fix \emph{all} governing characteristics; see \cite{Xu2009,barker16, guilleminot11,Lyckegaard2011,teferra18,laguerre2,Telley1992,Altendorf2014,kuehn08,Teferra2015,AFGK22}. The sites $\apdsites$ are generally chosen as the measured centroids of the grains. To model the shape of the grains, \cite{Lyckegaard2011} used Euclidean norms, i.e., $\apdmatrices=\Ecc$, while, following \cite{Alpers2015}, the heuristics of \cite{teferra18} use  covariance matrices. Further, the sizes~$\apdweights$ are typically estimated using the radius of the (Euclidean or ellipsoidal) ball whose volume coincides with the measured grain volume; see \cite{AFGK22} for an analysis in terms of discriminant analysis. 

As these heuristics abstain from any optimization whatsoever, they are computationally fast but generally result in a higher misclassification error. Hence, \cite{Alpers2015} suggested a different method that allows optimizing over $\Gamma$. The diagrams are specified through the dual variables of a linear program that addresses a constrained clustering problem; see \cref{sec:proxiess-GBPD} for details. As linear programming is often too time-consuming for the real-world instances of the desired size, we will, in particular, accelerate this method using data compression without any significant decrease of accuracy by employing \emph{coresets}. 

\subsubsection*{Contributions}
In the following three sections, we will 
\begin{itemize}
\item develop sparse variants of the GBPD-algorithm of \cite{Alpers2015}, called s-GBPD, which significantly improve running times,
\item extend the optimization to the set of all characteristic parameters within a data science framework which addresses $\APD$ directly leading to a \enquote{direct linear programming model} referred to as \DiM, and
\item state the results of a computational study on real-world data to indicate how the different algorithms perform in terms of computation times and accuracy.
\end{itemize} 
\cref{tab:problems2} provides an overview.

\begin{table}[h!]
    \centering
    \begin{tabular}{c|c|c|c|c}
    \toprule
  \multirow{2}*{Contributions} & \multicolumn{3}{c|}{Governing Characteristics} & \multirow{2}*{Section}\\ 
                      &              $\apdmatrices$ & $S$   & $\Gamma$ & \\ \midrule
\textbf{s-GBPD} & covariances & centroids & optimized & \ref{sec:proxiess-GBPD} \\[.5ex]
 \textbf{\DiM} &  optimized & optimized & optimized & \ref{sec:full-model}\\ \midrule
 \multicolumn{4}{l|}{Evaluation} & \ref{sec:evaluation_polycrystals} \\
         \bottomrule
    \end{tabular}
    \caption{Contributions: New algorithms, corresponding choice of characteristics and computational study. The entry \enquote{covariances} stands for the choices  $A_i=(\boldsymbol{\Sigma}_i)^{-1}$, where $\boldsymbol{\Sigma}_i$ is the covariance matrix of the $i$th grain.}
    \label{tab:problems2}
\end{table}

\section{Diagrams on sparse image supports}\label{sec:proxiess-GBPD}

In this section, we assume that the matrices $\apdmatrices$ and the sites $\apdsites$ are fixed beforehand (e.g., chosen as inverse covariance matrices and centroids, respectively), and only the sizes $\apdweights$ are subject to optimization. As, however, even $\APD(\apdmatrices,\apdsites)$ is computationally difficult, \cite{Alpers2015} introduced the \emph{linear programming} model GBPD. The computations can be viewed as a dual method that utilizes a relation between diagrams and weight-constrained anisotropic clustering. As shown by experimental computations in \cite{Alpers2015}, and later supported independently by \cite{sedivy16, sedivy17,Spettl2016}, GBPD is highly accurate and seems to capture the structure of polycrystals quite well. The computational cost, however, severely limits its practical use for 3D data sets of high resolution. Here we will show how to significantly accelerate the method, leading to a new class of algorithms, which we call sparse-GBPD or s-GBPD. In particular, we will show that the clustering can be performed at a rather low resolution, essentially without losing accuracy, and this will lead to significant speed-ups.

\subsection{Diagrams from weight-constrained clustering} 
\label{ssec:wcac}
Let us begin by briefly describing the relation between diagrams and clustering used in~\cite{Alpers2015}. For more details, we refer the reader to \cite{Brieden2017}, the handbook article \cite{GK17}, and the papers quoted therein.

As before, let $(X,Y)$ be a given grain scan, and let $\apdmatrices$ and $\apdsites$ be the covariance matrices and centroids derived from it. Further, we collect the measured cardinalities $\kappa_1,\ldots,\kappa_k,$ which, together with the resolution, approximate the grain volumes in the family $\Kappa=\{\kappa_1,\ldots,\kappa_k\}$. 

We will consider the more general case that each data point $x_j$ carries a weight $\omega_j \in (0, \infty)$. While each point in the original polycrystal carries the initial weight~$1$ we will use integer weights for compressed representations of multiple points by a single weighted point. (Let us mention in passing that, in other imaging applications, weights may, of course, be interpreted differently.) Collecting all weights in the family $\Omega=\{\omega_{1},\dots,\omega_{n}\}$, the weighted data set will then be denoted by $(X,\Omega)$. By convention, all parameters corresponding to the $j$th point will, again, be indexed by $j$. 

A \emph{clustering} $C$ of $(X,\Omega)$ is a vector
$$
C = (C_1,\ldots,C_k)=(\xi_{11},\dots, \xi_{1n}, \dots, \xi_{k1}, \dots, \xi_{kn})^\top
$$ 
whose components $\xi_{ij}$ specify the fraction of $x_j$ that is assigned to the \emph{$i$th cluster} $C_i=(\xi_{i1},\dots, \xi_{in})^\top$. More formally, the \emph{assignment conditions} are 
\begin{equation*} \label{eq:clustering_constraints}
\xi_{ij}  \ge 0 \quad \bigl(i \in[k],\ j \in [n]\bigr) \qquad \text{ and } \qquad \sum_{i=1}^{k} \xi_{ij}  = 1 \quad \big(j \in [n]\bigr).
\end{equation*}
Further, the \emph{weight} $\omega(C_i)$ of the cluster $C_i$ is given by
$$\omega(C_i)= \sum_{j=1}^{n} \xi_{ij}\omega_j.
$$
Note that the given grain scan $(X,Y)$ naturally defines the clustering $G=(G_1,\ldots,G_k)$ of the (discretized) grains. It is formally specified by the components $\xi^G_{ij}\in \{0,1\}$ with $\xi^G_{ij}=1$ if and only if $x_j \in G_i.$ As all points carry the weight $1$, cluster $G_i$ has weight
$\omega(G_i)= \kappa_i,$ $i\in[k].$ It is, of course, this \emph{grain scan clustering} that we want to represent by a suitable diagram $\apd=\apd(\apdmatrices,\apdsites, \Gamma)$. Accordingly, a clustering $C$ and a diagram $\apd=(P_1,\ldots,P_k)$ are called \emph{compatible} if, for all $i \in [k]$,
\begin{equation*}
    P_i \cap X \supset \text{supp}(C_i) = \{x_j \in X  \ : \ \xi_{ij} > 0 \}.
\end{equation*}
In general, we cannot expect the grain scan clustering to be compatible with any diagram hence we resort to the problem of determining a clustering $C^*$ that is \enquote{close} to $G$ and which admits a compatible diagram.

A natural and simple to model condition is to require that the clusters of $C^*$ and $G$ have similar weight, i.e., for given $\epsilon_{i} \in (0,\infty)$ the \emph{weight-constraints} 
\begin{equation*}
 (1-\epsilon_{i})\kappa_{i}  \le w(C_i) \le (1+\epsilon_{i})\kappa_{i} \qquad  \bigl(i \in [k]\bigr)
\end{equation*}
should hold. In practice, the error terms $\epsilon_{i}$ will reflect the available resolution of the image (or the known or expected accuracy of the measurements). As we allow fractional assignments, we may even set all $\epsilon_{i}$ to $0$ and regard the deviation on the weights inflicted by rounding to $0$-$1$-solutions as part of the approximation error. In the following, we will, therefore, enforce the strict constraints $w(C_i) =\kappa_{i}$ for all $i \in [k]$.
The set of all such \emph{weight-constrained clusterings} will be denoted by  $\Cs_{\Kappa}(k,X,\Omega)$. 
Further, the given shape parameters $\apdmatrices$ and $\apdsites$ are incorporated in the \emph{cost function}
\begin{equation*}
	\cost_\AK(X,C) = 
	\displaystyle \sum_{i=1}^{k}\sum_{j=1}^{n} \xi_{ij} \omega_j \norm{x_{j}-s_{i}}_{\apdmatrix_i}^{2}. 
\end{equation*}
The goal is now to find a weight-constrained clustering that minimizes this cost function, i.e., we are considering the problem \emph{weight-constrained anisotropic assigment} or, shorter,
\begin{quote}
    \wcaa:\\
    Given $k, X,\Omega,\Kappa$, find $C\in \Cs_\Kappa(k,X,\Omega)$ that minimizes $\cost_\AK(X,C)$.
\end{quote}
Since we can generally assume that $\Cs_\Kappa(k,X,\Omega) \ne \emptyset$, a finite minimum exists; it will be denoted by $\cost_\AK(X)$. Now, note that \wcaa can be formulated as the following linear program in the variables $\xi_{ij}$

\leqnomode
\begin{equation}
\begin{array}{ccrclcl}
	\tag{P}
	&\multicolumn{4}{c}{\displaystyle \min \quad \sum_{i=1}^{k} \sum_{j=1}^{n} \xi_{ij}\omega_{j}\norm{x_{j}-s_{i}}_{A_{i}}^{2}}&& \\[.6cm]
	&\qquad& 	\displaystyle \sum_{i=1}^{k} \xi_{ij}           & =   & 1          &\quad & \bigl(j \in [n]\bigr),\\[.4cm]
	&&	\displaystyle \sum_{j=1}^{n} \xi_{ij}\omega_{j} & =   & \kappa_{i} &\quad & \bigl(i \in [k]\bigr), \\[.4cm]
	&&	                                       \xi_{ij} & \ge & 0          &\quad & \bigl(i \in [k],\ j \in [n]\bigr).\\
\end{array}
\label{lp:primal}
\end{equation}
\leqnomode
Its dual 
\begin{equation}
\begin{array}{cccccccc}
 	\tag{D}
 	&           & \multicolumn{4}{c}{\displaystyle \max\quad \sum_{j=1}^{n} \omega_{j} \eta_{j} - \sum_{i=1}^{k} \kappa_{i}\apdweight_{i}} && \\[.6cm]
 	& \qquad & & \eta_{j} & \le & \norm{x_{j}-s_{i}}_{A_{i}}^{2} + \apdweight_{i}  &\quad & \bigl(i \in [k],\ j \in [n]\bigr),\\
 	\label{lp:dual}
 \end{array}
 \end{equation}
in the variables $\eta_{j}$ and $\apdweight_{i}$ encodes proximity conditions involving the diagram functions based on $\norm{x-s_{i}}_{A_{i}}^{2}$, and, in fact, we have the following theorem. 

\begin{theorem}[Special case of~{\cite[Thm.~1]{Brieden2017}}] \label{thm:wblsa:correspondence}
 	Let $(k, X,\Omega,\Kappa,\apdmatrices)$ be an instance of \wcaa, and $C\in \Cs_\Kappa(k,X,\Omega)$. Then $C$ is an optimizer of the linear program \labelcref{lp:primal} if and only if there exist size parameters $\apdweights$ such that $C$ and the diagram $\apd=\apd(\apdmatrices, \apdsites, \apdweights)$ are compatible.
\end{theorem}

Note that by linear programming duality, an optimal clustering $C^*$ results in an optimal diagram $\apd^*$. However, optimality is measured here with respect to the objective functions of~\cref{lp:primal} and~\cref{lp:dual}, and it is not clear in the first place how the cost of a clustering translates into the classification error which is to be minimized in $\APD$. (See \cite{FG22b} for related stability and instability results.) Hence, in effect, the clustering-based algorithms address~$\APD$ only indirectly.

\subsection{Data compression for accelerating the linear program} 
\label{sec:method:coresetandwblsa}

The linear program \cref{lp:primal} is based on the original ground set $X$ of the image (i.e., the points of $\nicefrac{1}{\rho}\cdot\Z^d$ in~$[0,1]^d$), resulting in $k\cdot (\rho +1)^d$ variables. If, for instance, in 3D, $\rho=9999$ and $k=50$, \cref{lp:primal} has already $5\cdot 10^{13}$ variables. Of course, as one would expect that $\xi_{ij}=0$ whenever $s_i$ is not among the first few closest sites to~$x_j$, preprocessing would allow reducing~$k$ to a manageable size. It is, however, $|X|$ which affects the number of variables the most. 

\subsubsection*{Image supports and coresets}

In the following, we will reduce the number of variables dramatically while retaining almost the initial accuracy. Our approach is based on the observation that the point set $(X,\Omega)$ is only a means for obtaining the desired diagram but not important in its own right. In fact, once we obtain the diagram $\apd$, we can easily generate an image at any desired resolution from it. Hence, we can, in principle, use \emph{any} non-empty finite subset $\supportX \subset [0,1]^d$ for computing diagrams through clustering -- as long as we  associate suitable weights to the points to accommodate the weight constraints. We will call such weighted sets $(\supportX,\supportOmega)$ \emph{image supports}. 

Naturally, some such sets $(\coresetX,\coresetOmega)$ will be better suited than others, and we will speak of \emph{coresets} if these sets can be derived from $(X,\Omega)$ in such a way that the relative deviation between $\cost_\AK(X)$ and $\cost_\AK(\coresetX)$ can always be suitably bounded in terms of a given precision. Here, we are interested in choosing the image supports $(\supportX,\supportOmega)$ as part of the  \enquote{experimental setup}, i.e., prior to scanning or taking characteristic measurements of the polycrystal. Thus, this setting requires the condition for \emph{all} sets of $k$ sites.

More precisely, let $(\coresetX,\coresetOmega)$ be a weighted data set, $\epsilon \in (0,\nicefrac{1}{2}]$, $\delta \in [1,\infty)$. Then, following \cite{FG22a}, $(\coresetX,\coresetOmega)$ is an \emph{$(\epsilon,\delta)$-coreset} for $(k,X,\Omega,\Kappa)$ if there exists a mapping $\extension:\Cs_\Kappa(k,\coresetX,\coresetOmega) \rightarrow \Cs_\Kappa(k,X,\Omega)$, called \emph{extension}, and real constants $\Delta^{+},\Delta^{-}$, referred to as \emph{$\Delta$-terms} or \emph{offsets}, with $0\le \Delta^{+} \le \delta \cdot \Delta^{-}$ such that the following two conditions hold for all sets of~$k$ sites and all clusterings $\coresetC\in  \Cs_\Kappa(k,\coresetX,\coresetOmega)$:
	\begingroup\leqnomode
	\begin{alignat}{3}
	(1-\epsilon)\cost_\AK(X,\extension(\coresetC)) &\le \cost_\AK(\coresetX,\coresetC) + \Delta^{+},  && \tag{a} \label{def:property_a} \\
	\cost_\AK(\coresetX) + \Delta^{-} &\le (1+\epsilon)\cost_\AK(X). \;  &&  \tag{b} \label{def:property_b}
	\end{alignat}
	\endgroup	
If $\delta=1$ we will simply speak of an \emph{$\epsilon$-coreset}.

Let us point out that other notions of coresets have been introduced for different problems; see e.g. \cite{Har-Peled2007,Feldman2013,FSS20,Fichtenberger2013,Sohler2018,Bachem2017a}. The definition used here is specifically customized for \wcacs and assignments. It guarantees that approximate solutions of the constrained  clustering problems on coresets can be turned into approximate solutions on the original data set; see \cite[Thm. 3.5]{FG22a}.

In order to practically use coresets for constrained clustering, they need to be small and the extension $g$ must be quickly computable. In our context, the extension does not play any role. In fact, clustering is not used in its own right but only as a method for computing diagrams via linear programming duality. Accordingly, we can determine the assignment of the original data points directly from the diagram computed from the coreset data and compare the obtained clustering $C$ with the grain scan clustering $G$ to quantify the classification error. 

In the following, we present two different coreset constructions which rely, in addition to $\apdmatrices$ and  $\apdsites$, only on the measurements $\kappa_1,\ldots,\kappa_k$. However, we will subsequently also consider an image support based on the grain scan clustering $G=(G_1,\ldots,G_k)$ to show how far the concept can be pushed. 

\subsubsection*{Pencil coresets} \label{sec:method:batching}

The following result shows that coresets exist and that they can be profitably utilized to accelerate GBPD computations. Let $\lambda^{-}(\apdmatrices)$ and $\lambda^{+}(\apdmatrices)$ denote the smallest and largest eigenvalue of all matrices in ~$\apdmatrices$, respectively.

\begin{theorem}[Special case of {\cite[Thm. 2.3]{FG22a}}]\label{th:smaller_coreset}
	For any instance of \wcaa and
	for any $\epsilon\in (0,\nicefrac{1}{2}]$ and 
	$$\delta \ge \frac{\lambdamax}{\lambdamin},$$
	there is an $(\epsilon,\delta)$-coreset $\coresetX$ of size 
	$$|\coresetX| \in \mathcal{O}\left(\frac{k^2}{\epsilon^{d+1}}\right).$$
\end{theorem}

The theorem shows that it suffices to use relatively small point sets for clustering whose \enquote{sparseness} is independent of the initial point set $X$ and depends only on $k$ and~$\epsilon$.  

Since we report on our implementation of pencil coresets for different error parameters in \cref{sec:evaluation_polycrystals}, we will briefly indicate their structure. The construction is based on a design by \cite{Har-Peled2007} for unconstrained clustering; its generalization to \wcacs and improved bound is due to \cite{FG22a}. 

The pencil coresets \enquote{live} on pencils of lines issuing from a certain set of points constructed in a first step. For technical simplicity, we will focus on the Euclidean case and, as this is the relevant situation here, assume that these points coincide with given sites $s_1,\ldots,s_k$. The lines of the $i$th pencil are then chosen in such a way that their points of intersection with the unit sphere centered at the site $s_i$ are \enquote{evenly distributed}. Note that this is closely related to coverings of the unit sphere by congruent spherical caps. See, e.g., \cite{M02} for general results on \emph{$\epsilon$-nets} and \cite{F72} for sparse coverings of the Euclidean $2$-sphere by specific numbers of congruent spherical caps. 

Then, we assign each point of the data set to its closest site, breaking ties arbitrarily, and project it onto a line of the corresponding pencil closest to it, again breaking ties arbitrarily. After this projection process, all points lie on lines, and in the remaining step, each line is treated independently. Note that this is possible due to the additivity of the cost function for the data points. 

Each such line is now partitioned into batches. Finally, all points in one batch are replaced by the batch's centroid and assigned the sum of weights of all batch points as its weight. 

If the lengths of the batches are chosen appropriately, this process results in the desired $(\epsilon,\delta)$-proxy. See \cite{FG22a} for the details of the proof.

\subsubsection*{Resolution coresets} \label{sec:method:resolution}

While the construction in the previous section works for general data sets $(X,\Omega)$, the pencils ignore the grid structure of the given grain scan. In particular, the \enquote{density} of the image support varies with the distance to the sites.
For voxel-based images, it seems much more natural and, actually, interesting in its own right, to consider image supports that pertain to a (coarser) grid structure and which can still be interpreted in terms of the image resolution.

In the following, we will briefly introduce the concept. For technical simplicity, we will assume that the voxels are boxes whose edge lengths are powers of $\nicefrac{1}{2}$. Our exposition follows \cite{FGK22}.

For $i \in [d]$, and $\rho= (\rho_1,\dots, \rho_d) \in \N^d$ we consider the uniform partitioning of the interval $[0,1]$ on the $i$th coordinate axis into $2^{\rho_i}$ intervals of equal length. Then each voxel is the Cartesian product of one of the intervals on each axis. 
More formally, with
\begin{equation*}
x_{j_1,\ldots,j_d}^{\rho}= \left(
\begin{array}{c}
\displaystyle \frac{1}{2^{\rho_1+1}} + \frac{j_1}{2^{\rho_1}} \\
\vdots\\
\displaystyle\frac{1}{2^{\rho_d+1}} + \frac{j_d}{2^{\rho_d}}\\
\end{array}
\right),
\end{equation*}
the voxels are obtained as
\begin{equation*}
x_{j_1,\ldots,j_d}^\rho
+ \frac{1}{2^{\rho_1+1}} [-1,1] \times \ldots \times \frac{1}{2^{\rho_d+1}} [-1,1]
\end{equation*}
for all choices of $j_1,\ldots,j_d$ with
\begin{equation*}
j_1\in \{0\} \cup [2^{\rho_1}-1], \ldots, j_1\in \{0\} \cup [2^{\rho_d}-1].
\end{equation*}
We collect the points $x_{j_1,\ldots,j_d}^{\rho}$ in the data set $\gridX^{\rho}$. Note that the points are centroids of the voxels and constitute the original data set, i.e., $\gridX^{\rho} = X$. Of course, $n=|X^{\rho}|= 2^{\rho_1}\cdot \ldots \cdot 2^{\rho_d}$. As before, we will refer to the vector $(\rho_1, \ldots, \rho_d)$ as the \emph{resolution} of the image. Note that, if all $\rho_i$ coincide with some given $\bar\rho \in \N$ then 
\begin{equation*}
    \gridX^{\rho}= \frac{1}{2^{\bar\rho+1}} \1 + \left(\left(\frac{1}{2^{\bar\rho}}\Z^d\right) \cap [0,1)^d\right). 
\end{equation*}
Naturally, a lower resolution is obtained by decreasing some or all $\rho_i$, resulting in sparser image supports. The question is how small we can choose a resolution $\tau \in \N^d$ with $\tau_i \le \rho_i$ such that $\gridX^{\tau}$ is a coreset -- together with appropriate weights.

As it turns out, it suffices to consider a number of points on each axis which depends only on $k$ and $\epsilon$.

\begin{theorem}[ {\cite[Thm.~1]{FGK22}}]\label{th:resolution-coreset}
	Let $\rho \in \N^d$, $(k,\gridX^\rho,\Kappa)$ be an instance of \wcac and
  $\epsilon\in (0,\nicefrac{1}{2}]$. Then, there is a $\bar\tau \in \N $ with
	\begin{equation*}
	    \bar\tau \le \frac{2^{\nicefrac{8}{3}}k}{\epsilon^{\nicefrac{2}{3}}},
	\end{equation*}
	such that, for $\tau =(\bar\tau,\dots, \bar\tau) \in \N^d$, the set $X^{\tau}$ together with appropriate weights is an $\epsilon$-coreset of size 
	$$
	|X^{\tau}| \le \left(\frac{2^{\nicefrac{8}{3}}k}{\epsilon^{\nicefrac{2}{3}}}\right)^d.
	$$
\end{theorem}

\Cref{th:resolution-coreset} allows us to reduce the resolution of the scan to a size that is independent of the original resolution and only depends on the number of grains and the desired accuracy (measured in terms of the cost of the produced clustering). Hence, we can use grain scans of relatively low resolution, solve \cref{lp:primal} for the correspondingly small data set, and then obtain a diagram representation via \cref{lp:dual}. Recall that we can, ultimately, use the diagram to produce a grain image of any (arbitrarily high) desired resolution.

Let us briefly address the difference between the bounds in \cref{th:smaller_coreset} and \cref{th:resolution-coreset} to place the results into perspective. For example, when $k=50$ and $\epsilon= 0.01$, not taking constants (hidden, in \cref{th:smaller_coreset}, in the Landau notation) into account, the coreset size guaranteed by \cref{th:smaller_coreset,th:resolution-coreset} in 3D is $2.5\cdot 10^{11}$ and $1.25\cdot 10^{9}$, respectively. Note, however, that \cref{th:smaller_coreset,th:resolution-coreset} give \emph{worst-case upper} bounds which are not known (and not even conjectured) to be best possible. As we will see in \cref{sec:evaluation_polycrystals}, the image supports can indeed be reduced much more radically in our practical computations without sacrificing much of the accuracy. 

Let us further remark that, while the coreset property bounds the error in clustering costs it does not directly provide a theoretical bound for the classification error of the corresponding diagram for  $\APD(\apdmatrices,\apdsites)$. The practical results in \cref{sec:evaluation_polycrystals} indicate, however, that the \enquote{indirect} quality guarantee which comes with the coreset property also leads to a favorable classification behavior.

\subsubsection*{Grain scan-based interior removal}

The previously constructed pencil and resolution coresets make only use of certain characteristics, which, of course, can be derived from a given grain scan $(X,Y),$ but which, in practice, are often also available through direct measurements not requiring full knowledge of the grain scan clustering $G=(G_1,\ldots,G_k)$. While this is an important feature for practical applications it is, of course, natural to ask how much is actually lost by not utilizing the full grain scan information. More precisely, using $G$, how much further can the size of image supports be reduced while still constituting coresets, i.e., while still facilitating quite accurate clusterings in terms of their costs. Note that this is not just a theoretical issue but also highly relevant whenever given grain scan representations need to be turned into diagram representations to facilitate image post-processing and analysis tasks.

To set the stage, let us mention first what happens if not just $(X,Y)$ but, in addition, a best fitting diagram $\apd^*$ is available. Then, by \cite{FG22a}, $\apd^*$  is determined already by at most $\mathcal{O}(k^2)$ points. In fact, points in the scan closest to the cell boundaries suffice, where closeness is measured with respect to the local ellipsoidal norms. Thus, if we knew an optimal diagram, we could identify such points and solve \wcaa using only $\mathcal{O}(k^2)$ data points, i.e., via the linear program \cref{lp:primal} in only $\mathcal{O}(k^3)$ variables.

Of course, best fitting diagrams are generally not known beforehand, and it is precisely our task to compute them. Now, suppose that, under ideal conditions, the grain growth process indeed results in a diagram. Then we may assume that the given grain scan clustering~$G$ is close, up to minor imperfections, to an optimal clustering $C^*.$ 

Before continuing, let us mention that such an assumption is in line with common stochastic grain growth models based on the normal distribution $\mathcal{N}(\mu_i,\Sigma_i)$ for the points in the $i$th grain. Here, the shape parameters are chosen according to maximum likelihood resulting in $s_i=\mu_i$ and $A_i=\Sigma_i^{-1}$. See \cite{AFGK22} for more information. 

Under such a \enquote{closeness} assumption it seems reasonable to discard all points that are sufficiently far away from the \enquote{grain boundaries} of $G$. More formally, we choose a metric $\mathbbm{d}$ and parameters $\delta_i \in (0,\infty)$, and define the \emph{$\delta_i$-interior} of $G_i$ by 
\begin{equation*}
	G_{i}(\delta_i) = \{x_j \in G_{i} :
 \mathbbm{d}(x_{j}, x_{\ell}) \ge \delta_i \, \, \forall \ell \in [k]: y_{\ell} \ne i  \}.
\end{equation*}
First note that, for the Euclidean norm, this is a discrete variant of the well-known concept of an \emph{inner parallel body}. In the grid based situation of the previous paragraph it seems more natural to use a variant of the $1$-norm normalized by the different resolutions on the axes, i.e.,
for $z=(\eta_1,\ldots,\eta_d)^\top$
\begin{equation*}
    \norm{z}_{1}= \sum_{i=1}^d 2^{\rho_i}\eta_i.
\end{equation*}
The corresponding metric then counts the \enquote{grid distance} which is the length of a shortest path in the underlying grid graph with edge weights $1$. In \cref{sec:evaluation_polycrystals}, we will use this metric and refer to it as \emph{grid graph distance}.
Note, that (for any standard metric)~$G_{i}(\delta_i)$ can be computed efficiently. See \cite{FG22a} for a more comprehensive account on image supports obtained this way.

The parameters $\delta_i$ control the size of the interior, and a larger $\delta_i$ will result in a smaller $\delta_i$-interior and, hence, in a larger image support but potentially also higher accuracy. 
The specific choices of the $\delta_i$ reflect the initial confidence that $G$ is already close to an optimal clustering. The deviation from $G$ of the grain scan induced by the computed diagram $\apd$ may be viewed as an indicator of the state of the process, i.e., of how far it is from an energy-minimal state which is indeed characterized by the underlying model.

For the computations in \cref{sec:evaluation_polycrystals}, we employ image supports obtained by coreset techniques combined with interior removal. Any specific algorithm of the so obtained class of methods will be referred to as \emph{sparse-GBPD} or, short, \emph{s-GBPD}. As we will show, suitable such choices lead to substantial accelerations of GBPD computations on real-world data. 

\section{A direct model for $\APD$}
\label{sec:full-model}

The algorithms of the previous section are based on estimates for the matrices and sites and address $\APD(\apdmatrices, \apdsites)$ indirectly via clustering. We will now propose an optimization model which solves $\APD$ directly. As it optimizes over all three characteristic parameters $\apdmatrices$, $\apdsites$, and $\apdweights$, the larger set of diagrams to choose from will allow for a better fit of~$(X,Y)$. The higher accuracy comes, however, at the price of computational costs which are substantially higher than those of the methods described in \cref{sec:proxiess-GBPD}. Nonetheless, combining the model with our previous point reduction techniques enables at least its computation. 

Hence, this section has a somewhat different focus: On the one hand, it provides a \enquote{benchmark method} for a better assessment of the quality of fit produced by other methods for appropriately chosen test cases. On the other hand, it may be helpful for \enquote{quantifying the physical state} of a prematurely terminated growth process in terms of the achieved fit. As expected and exemplified in \cref{sec:evaluation_polycrystals}, the new algorithm will only be practically useful if computation time is not the main concern.

The technique that we will now introduce is based on the paradigm of \emph{support vector machines} from data analysis and requires knowledge of (a relevant part of) the grain scan~$(X,Y)$. See, e.g.,~\cite{SS02} for additional background information on machine learning.

Let, again $G=(G_1,\ldots,G_k)$ be the grain scan clustering. If there exists a diagram~$\apd(\apdmatrices, \apdsites, \apdweights)$,  compatible with $G$, its defining functions $h_i(x) = \norm{x-s_{i}}_{\apdmatrix_i}^2 + \apdweight_{i}$ will separate the points of the clusters, i.e., for $j\in [n]$ and $i,\ell\in [k]$ with $\ell\ne i$,
\begin{equation*}
    h_i(x_j) \le  h_\ell (x_j) \qquad \text{for each $x_j \in G_i$}.
\end{equation*}
Since, in general, real-world scans will not admit a perfectly fitting diagram, it is natural to relax these conditions and resort to the optimization problem
\begin{equation*}
\begin{array}{ccrclcl}
	&\multicolumn{4}{c}{\displaystyle \min \quad \sum_{j=1}^{n}\zeta_j}&& \\[.6cm]
	&&	\displaystyle h_i (x_j) & \le  & h_\ell (x_j) + \zeta_j &\quad & \bigl(j\in [n], \, i,\ell \in [k], \, \ell\ne i, \, x_j\in G_i\bigr) \\[.4cm]
	&&	              \zeta_{j} & \ge  & 0                      &\quad & \bigl(j \in [n]\bigr).\\
\end{array}
\end{equation*}
Here, the non-negative error terms $\zeta_j$ allow for potential misclassification, and the objective minimizes the accumulated errors. As, however, we can multiply all matrices $A_i$ and sizes $\gamma_i$ by a common positive factor without changing the diagram, the infimum of the objective function will always be $0$, independently of the fit of the diagram. There are various possibilities to break this invariance. Our approach, reflects the assumption that points near the grain centroids are already correctly classified by the grain scans and classification errors occur only near the boundaries. This will be modeled in terms of the $\delta_i$-interior for suitably chosen $\delta_i \in (0,\infty)$ introduced in the previous section, resulting in the optimization problem
\begin{equation*}
\begin{array}{ccrclcl}
	&\multicolumn{4}{c}{\displaystyle \min \quad \sum_{j=1}^{n}\zeta_j}&& \\[.6cm]
	&&	\displaystyle h_i (x_j) & \le  & h_\ell (x_j) -1 &\quad & \bigl(j\in [n], \, i,\ell \in [k], \, \ell\ne i, \, x_j\in G_i(\delta_i)\bigr) \\[.4cm]
	&&	\displaystyle h_i (x_j) & \le  & h_\ell (x_j) + \zeta_j &\quad & \bigl(j\in [n], \, i,\ell \in [k], \, \ell\ne i, \, x_j\in G_i\setminus G_i(\delta_i)\bigr) \\[.4cm]
	&&	              \zeta_{j} & \ge  & 0          &\quad & \bigl(j \in [n], \, x_j\in G_i\setminus G_i(\delta_i)\bigr).\\
\end{array}
\end{equation*}
In this model, all points in the $\delta_i$-interior of $G_i$ are separated strictly. Note that the specific choice of $1$ for the \enquote{margin} is arbitrary as long as it is positive.

Unfortunately, this optimization problem has nonlinear constraints as the optimization extends over the entries in $A_i$, $s_i$ and $\gamma_i$ and
\begin{equation*}
	h_i(x_j) = \norm{x_j-s_{i}}_{\apdmatrix_i}^2 + \apdweight_{i} 
	= x_j^TA_ix_j-2s_i^TA_ix_j+s_i^TA_is_i+\gamma_i,
\end{equation*}
where only the points $x_j$ are given. Using the notation
\begin{equation*}
    x_j= \left(
    \begin{array}{c}
         (x_j)_1  \\
         \vdots \\
         (x_j)_d \\
    \end{array}
    \right), \qquad
    A_i= \left(
    \begin{array}{ccc}
         (A_i)_{1,1} & \dots & (A_i)_{1,d}  \\
         \vdots && \vdots\\
         (A_i)_{d,1} & \dots & (A_i)_{d,d}  \\
    \end{array}
    \right)
\end{equation*}
and setting
\begin{equation*}
    a_{i} = -2\apdmatrix_{i}s_{i}, \qquad \alpha_{i} = s^{\top}_{i}\apdmatrix_{i}s_{i}+\apdweight_{i}
\end{equation*}
we can encode all parameters in the vector
\begin{equation*}
	\fa_{i} = \Bigl(\alpha_{i},a_{i}^{\top},(\apdmatrix_{i})_{1,1},2(\apdmatrix_{i})_{1,2},\dots,2(\apdmatrix_{i})_{d-1,d},(\apdmatrix_{i})_{d,d}\Bigr)^{\top}\in \R^{1+d+\frac{1}{2}d(d+1)}.
\end{equation*}
Note that $(A_i,s_i,\gamma_i)$ can easily be recovered from $\fa_{i}$. In fact, the entries of $A_i$ can be read off directly, then $s_{i} = -\nicefrac{1}{2}\apdmatrix_{i}^{-1}a_{i}$ and, finally, $\apdweight_{i} = \alpha_i-s^{\top}_{i}\apdmatrix_{i}s_{i}$.
If we now introduce, for each point $x_j$, its extended coefficient vector 
\begin{equation*}
	\fx_j = \Bigl(1,(x_j)_{1},\dots,(x_j)_{d},(x_j)_{1}^{2},(x_j)_{1}(x_j)_{2},\dots,(x_j)_{1}(x_j)_{d},\dots, (x_j)_{d-1}(x_j)_{d},(x_j)_{d}^{2}\Bigr)^{\top},
\end{equation*}
we can write the separation conditions in terms of the linear expression 
\begin{equation*}
  h_{i}(x_j) = \fa_{i}^{\top}\fx_j. 
\end{equation*}
Accordingly, we will address $\APD$ via the \emph{direct linear programming model} 
\leqnomode
\begin{equation*}
\begin{array}{ccrclcl}
\tag{\DiM}
	&\multicolumn{4}{c}{\displaystyle \min \quad \sum_{j=1}^{n}\zeta_j}&& \\[.6cm]
	&&	\displaystyle \fa_{i}^{\top}\fx_j - \fa_{\ell}^{\top}\fx_j +1  & \le  & 0 &\quad & \bigl(j\in [n], \, i,\ell \in [k], \, \ell\ne i, \bigr. \\
		&&	&  &&\quad & \bigl. \,\, x_j\in G_i(\delta_i)\bigr) \\[.2cm]
	&&	\displaystyle \displaystyle \fa_{i}^{\top}\fx_j - \fa_{\ell}^{\top}\fx_j - \zeta_j & \le  & 0 &\quad & \bigl(j\in [n], \, i,\ell \in [k], \, \ell\ne i, \bigr. \\
		&&	&  &  & & \bigl. \,\, x_j\in G_i\setminus G_i(\delta_i)\bigr) \\[.2cm]
	&&	              \zeta_{j} & \ge  & 0          &\quad & \bigl(j \in [n], \, x_j\in G_i\setminus G_i(\delta_i)\bigr).\\
\end{array}
\label{lp:polynomial_boundaries_II}
\end{equation*}
In the following we will also use the acronym \DiM for any algorithm that solves (\DiM). 

There is one additional issue to address. In order to derive diagram parameters from a solution $(\fa_1^*,\ldots,\fa_k^*)$ of \cref{lp:polynomial_boundaries_II}, we need to guarantee that the obtained matrices $A_1^*,\ldots,A_k^*$ are symmetric and positive definite. Otherwise, the computations may not result in an anisotropic power diagram. While the symmetry of the matrices is built into the model, there is, however, no precaution taken that the $A_i^*$ produced via~\DiM are indeed positive definite. 

Suppose now that not all of the resulting matrices $A_i$ are already positive definite. Then we simply switch over to $\bar{A}_i=A_i+\beta E_d$ for all $i\in [k]$ and some suitable $\beta \in (0,\infty)$ so that all $\bar{A}_i$ become positive definite. Note that the constraints in \cref{lp:polynomial_boundaries_II} do not change. In effect, we simply add and subtract the term $\beta x_j^{\top}x_j$ in the first two groups of constraints. Hence, we can use the matrices $\bar{A}_i$ for the subsequent computation of the sites and the sizes and set
\begin{equation*}
     \bar{s}_{i} = -\nicefrac{1}{2}\bar{\apdmatrix}_{i}^{-1}a_{i}, \qquad  \bar{\apdweight}_{i} = \alpha_i-s^{\top}_{i}\bar{\apdmatrix_{i}}\bar{s}_{i}.
\end{equation*}
Then $(\bar{A}_i,\bar{s}_i,\bar{\gamma}_i)$ ($i\in [k]$) yield, indeed, a desired diagram $\apd(\bar{\apdmatrices}, \bar{\apdsites}, \bar{\apdweights})$ with the same objective function value in \cref{lp:polynomial_boundaries_II}.

Let us close by pointing out that \cref{lp:polynomial_boundaries_II} can be combined with the image supports described in \cref{sec:method:coresetandwblsa}. Indeed, for practical computations one can largely thin out~$X$ first and, in particular, require strict separation only for data points in certain rings $G_i(\delta^-_i)\setminus G_i(\delta^+_i)$ with $\delta^-_i < \delta^+_i$. In \cref{sec:evaluation_polycrystals} we will employ such accelerating techniques but yet see that \cref{lp:polynomial_boundaries_II} is significantly more time consuming than \cref{lp:primal}.

\section{Computational results}
\label{sec:evaluation_polycrystals}
We will now evaluate the introduced techniques with respect to computation time and accuracy using a 3D real-world grain scan from \cite{Lyckegaard2011a}. After briefly describing the computational setup, we will first assess the achievable degree of sparseness of the image supports for s-GBPD. Then, we will use the so obtained \enquote{practical gauging} of the \enquote{sparseness parameters} to compare s-GBPD and \DiM for a range of different quality criteria. This comparison will illustrate the degree of fit produced by the former at a much lower computational cost than the latter.

\subsection{Computational setup}
Our computational study is based on data provided by \cite{Lyckegaard2011a}; see also the description in \cite{Alpers2015}. 
It was obtained by a synchrotron mico-tomograph experiment conducted on a metastable beta titanium alloy (Ti $\beta$21S). The material was scanned with a volume of 240\textmu m $\times$ 240\textmu m $\times$ 420\textmu m at a resolution resulting in a voxel edge length of 0.7~\textmu m. It comprises $591$ grains, $211$ are interior to the polycrystal, i.e., they are entirely surrounded by other grains.

\begin{figure}[h!]
	\centering
	\includegraphics[width=.7\linewidth]{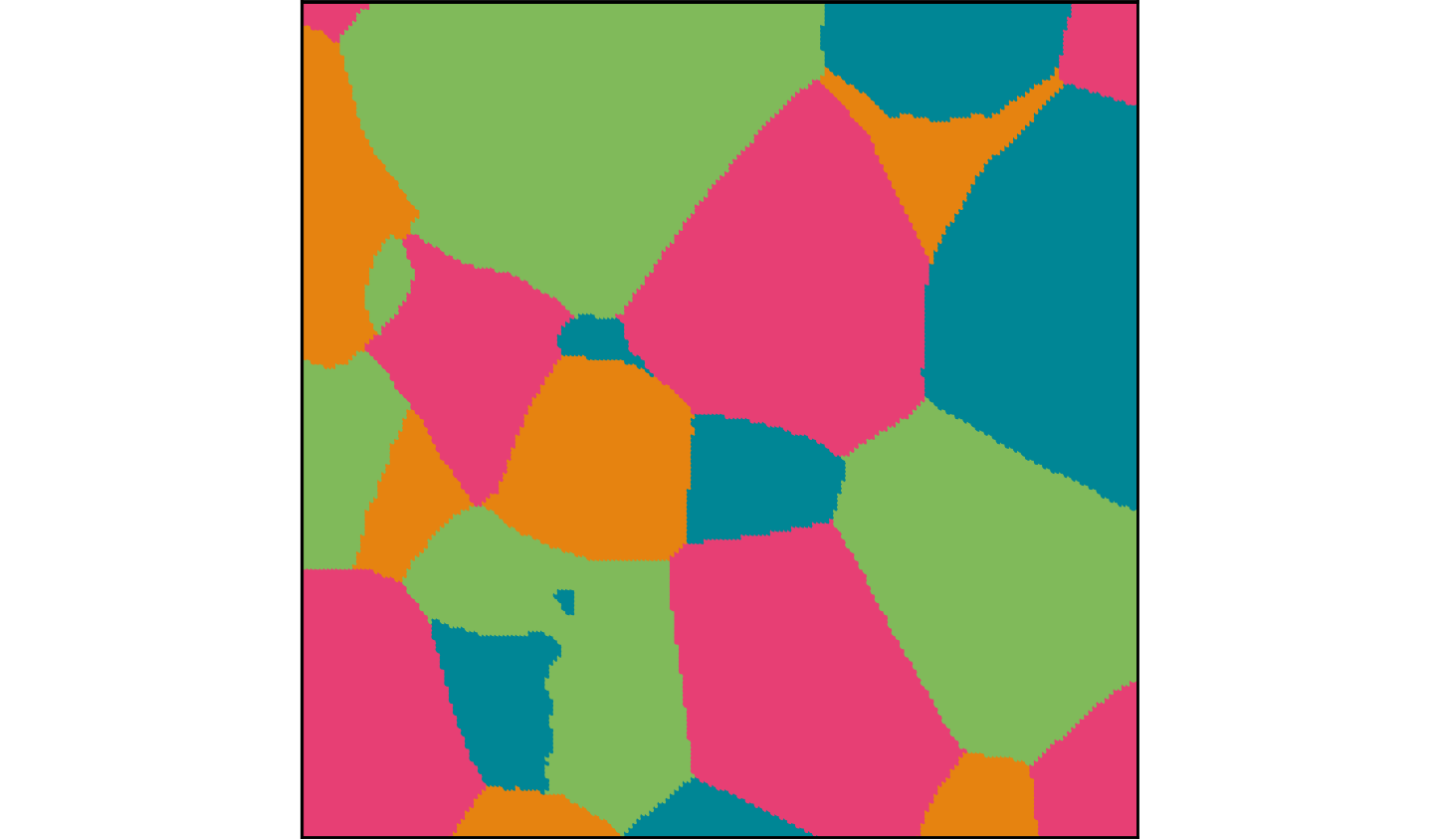}
	\caption{Slice no. 50 of the data set. Neighboring grains are colored in different colors.}
	\label{fig:grainscan-2d}
\end{figure}

We hence have a grain scan $(X,Y)$ with
\begin{equation*}
    n=|X|= 339 \times 339 \times 599 = 68,837,679 \qquad \textnormal{and}\qquad k=591.
\end{equation*}
In the following we refrain from normalizing $X$ and report all results for the original physical dimension stated above. All computations are carried out in 3D. The images depict 2D slices parallel to the plane of the first two coordinates. One such example is shown in \cref{fig:grainscan-2d}. 

Let us point out that the larger blue patch in the left lower corner and the close-by tiny blue patch surrounded by the green area actually belong to the same grain. While this can be seen in the 3D image, its 50th slice depicts two components. So, let us emphasize again that all subsequent 2D images are for illustration purposes only while all reported results have been fully computed in 3D. This contrasts most computations in \cite{Alpers2015} and, as a matter of fact, 3D problems of the size considered here are far beyond the capability of the original GBPD algorithm introduced there. 

The reported running times were measured on a virtual machine with ten virtual CPUs corresponding to five Intel(R) Xeon(R) Gold 6148 CPU @ 2.40GHz CPUs and 45GB of memory (i.e., one CPU processes two virtual CPUs). For better comparison, we solved all linear programs using the dual simplex algorithm on only one virtual CPU. Let us point out, however, that the sparse models, which we show to be sufficient, also run on standard machines with less than 16GB of memory.

\subsection{s-GBPD: Sparseness vs. fit}
The coreset constructions in Theorems \ref{th:smaller_coreset} and \ref{th:resolution-coreset} allow us to work with sample sizes that only depend on $k$ and $\epsilon$ without much affecting the clustering cost. Grain-based interior removal can further reduce the size of the image support. We will now study experimentally the effect of sparseness on the quality of the generated diagram. 

\subsubsection*{Accuracy and relative cluster weight error}
In the following, we measure the quality of fit. More precisely, let $\apd=(P_1,\ldots,P_k)$ be the diagram obtained by the algorithm and set  
\begin{equation*}
  \hat{C}_i= \bigl\{x_j\in X: x_j \in \interior(P_i)\bigr\} \quad \bigl(i\in [k]\bigl), \qquad \hat{C}=(\hat{C}_1,\ldots,\hat{C}_k).
\end{equation*}
Let us point out that $\hat{C}$ is not a clustering of $X$ but only of $X\cap \bigcup_{i=1}^k \interior(P_i)$ (and this is why we use the notation $\hat{C}$ rather than $C$). We define the \emph{relative fit} or \emph{accuracy} $\Phi_G$ of $\hat{C}$ as
\begin{equation*}
  \Phi_G(\hat{C})= \frac{1}{n} \left| \bigcup_{i=1}^k G_i\cap \hat{C}_i \right|.
\end{equation*} 
Note that the points on the cell boundaries are not assigned to any clustering and, hence, count as misclassified. 

Also, we will compute the \emph{relative cluster weight error} $\Psi_G$ defined by
\begin{equation*}
	\Psi_G(\hat{C}) = \frac{1}{n}\sum_{i=1}^k|\kappa_i - \omega(\hat{C}_i)|.
\end{equation*}

\subsubsection*{Pencil coresets} 

In this setup, the sites $\apdsites$ and the matrices $\apdmatrices$ are given and we consider pencils issuing from the centroids of the $G_i$ whose lines are distributed with respect to the ellipsoids $s_i+\{x:\|x\|_{A_i}\le 1\}$.

\begin{figure}[h!]
	\centering
	\includegraphics[width=1\linewidth]{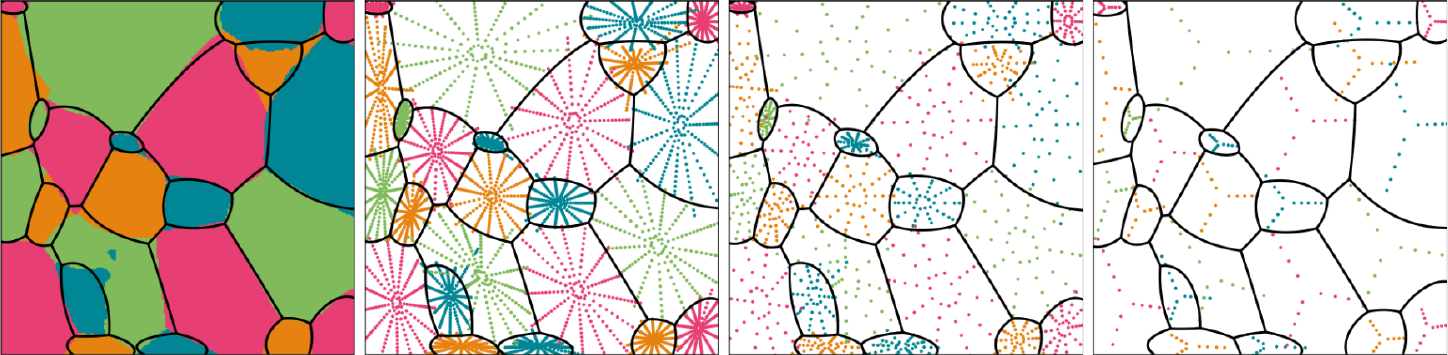}
	\caption{GBPDs from pencil coresets. From left to right: Computed diagrams (cell boundaries shown in black) based on decreasing image supports (the non-white and non-black colored pixels in the respective images). The diagram in the left image uses the full image support (the full available data), while the others use image supports living on pencils of decreasing density.}
	\label{fig:grainscan-2d-batchingcoreset}
\end{figure}

While all computations are fully 3D, but this is hard to visualize,  \cref{fig:grainscan-2d-batchingcoreset} illustrates the construction in 2D with rays emanating from the centroids of the depicted 2D grains. 

\begin{figure}[h!]
	\centering
	\includegraphics[width=1\linewidth]{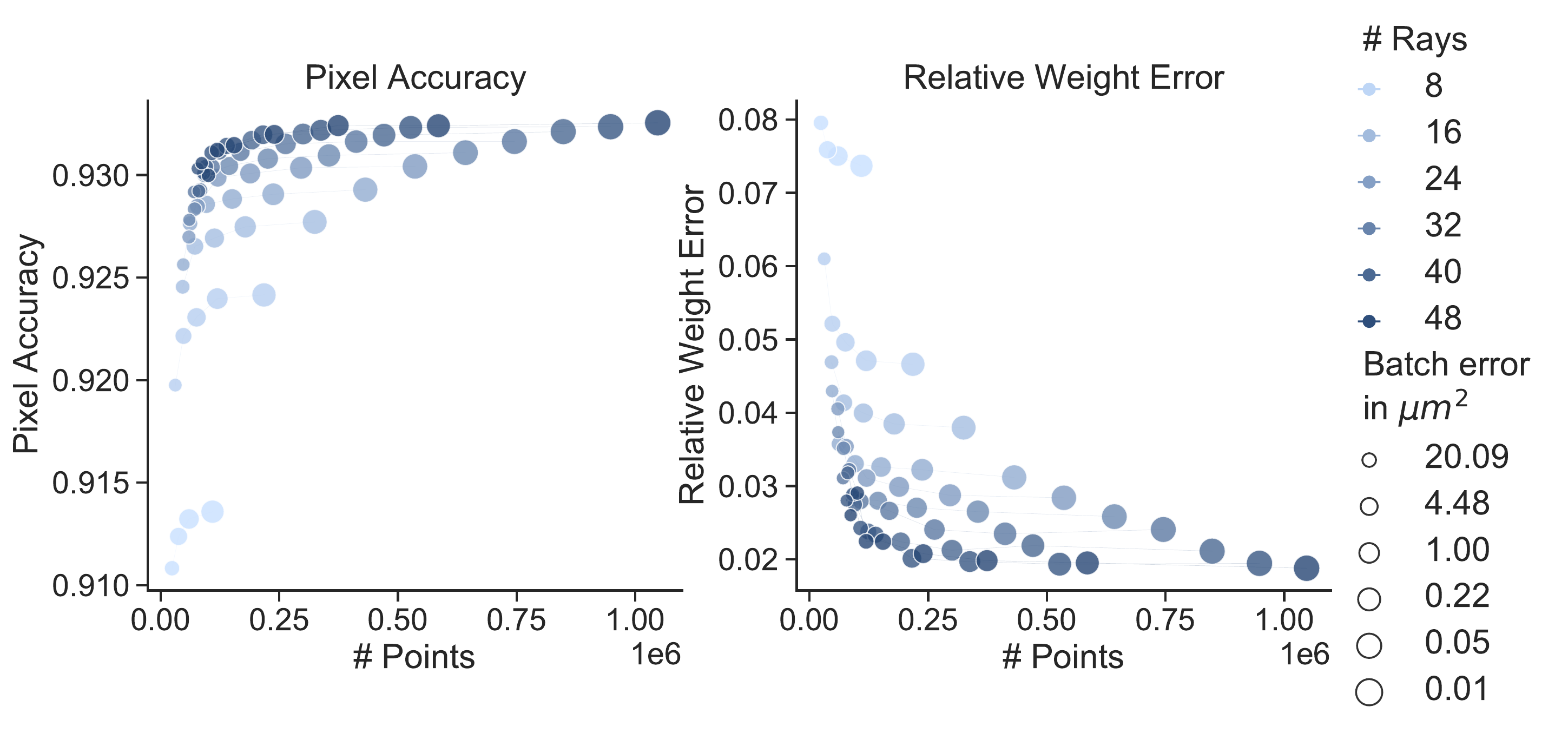}
	\caption{Accuracy (left) and relative cluster weight error (right) for different cardinalities of image supports. The numbers of rays and the batch errors are signified by color and size. Also, the behavior for a fixed number of rays is indicated by connecting the corresponding points.}
	\label{fig:3dgrainscanbatchingfastwblsa}
\end{figure}

For about 100,000 points, we obtain already an accuracy of $0.93$. Enlarging the image support improves $\Phi_G$ only in the third decimal digit. The relative cluster weight error $\Psi_G$ benefits slightly more from increasing the density, from 0.025 for 100,000 to 0.02 for 200,000 points. \cref{fig:3dgrainscanbatchingfastwblsa} depicts the dependency of the accuracy and the relative cluster weight error on the density of the sample support. While in \cref{th:smaller_coreset}, both, the number of rays and the batch error are functions of $k$ and $\epsilon$, we control these numbers in practice separately.

More precisely, first the number of rays are selected and each point of $X$ is projected on the appropriate ray, as described in \cref{sec:method:batching}. Then the points of each ray are partitioned into batches $B$ and replaced by the batches' centroids $c$. We do this dynamically so as to control the batch errors $\sum_{x \in B}\norm{x-c}_{A}^{2}$. Accordingly, \cref{fig:3dgrainscanbatchingfastwblsa} also shows the influence of the batch error and the number of rays (indicated by the size and color of the points, respectively) on the pixel accuracy and the relative weight error.

In general, the computational results suggest that the number of rays may be more important than the batch error. This observation supports the intuition that the shape information about the grains is better captured by additional rays than by more points on each ray.

\subsubsection*{Resolution coresets}
\begin{figure}
	\centering
	\includegraphics[width=1\linewidth]{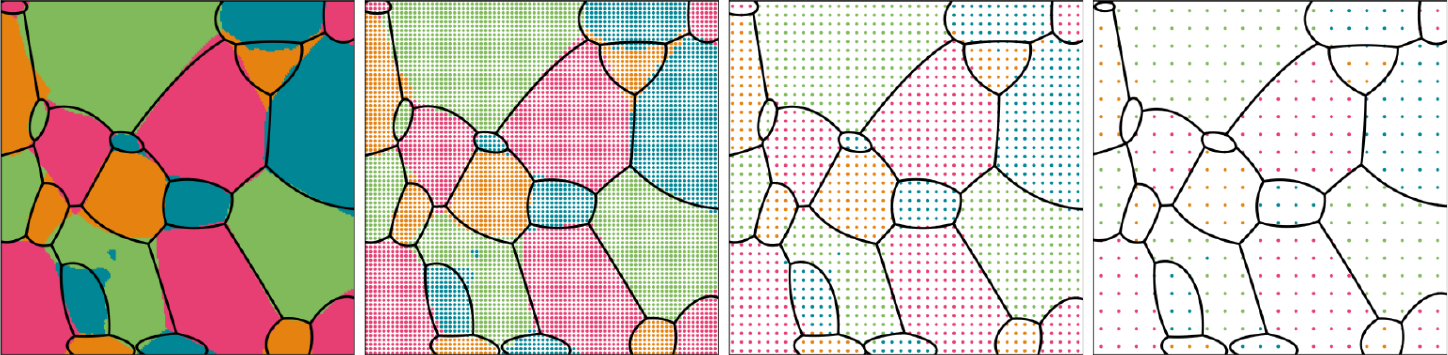}
	\caption{GBPDs from resolution coresets. From left to right: Computed diagrams (cell boundaries shown in black) based on decreasing image supports (the non-white and non-black colored pixels in the respective images). The diagram in the left image uses the full image support (the full available data), while the others use image supports living on resolution coresets of decreasing density.}
	\label{fig:grainscan-2d-resolutioncoreset}
\end{figure}

The results of a similar evaluation for the resolution coresets are illustrated in \cref{fig:grainscan-2d-resolutioncoreset}.

\begin{figure}[h!]
	\centering
	\includegraphics[width=1\linewidth]{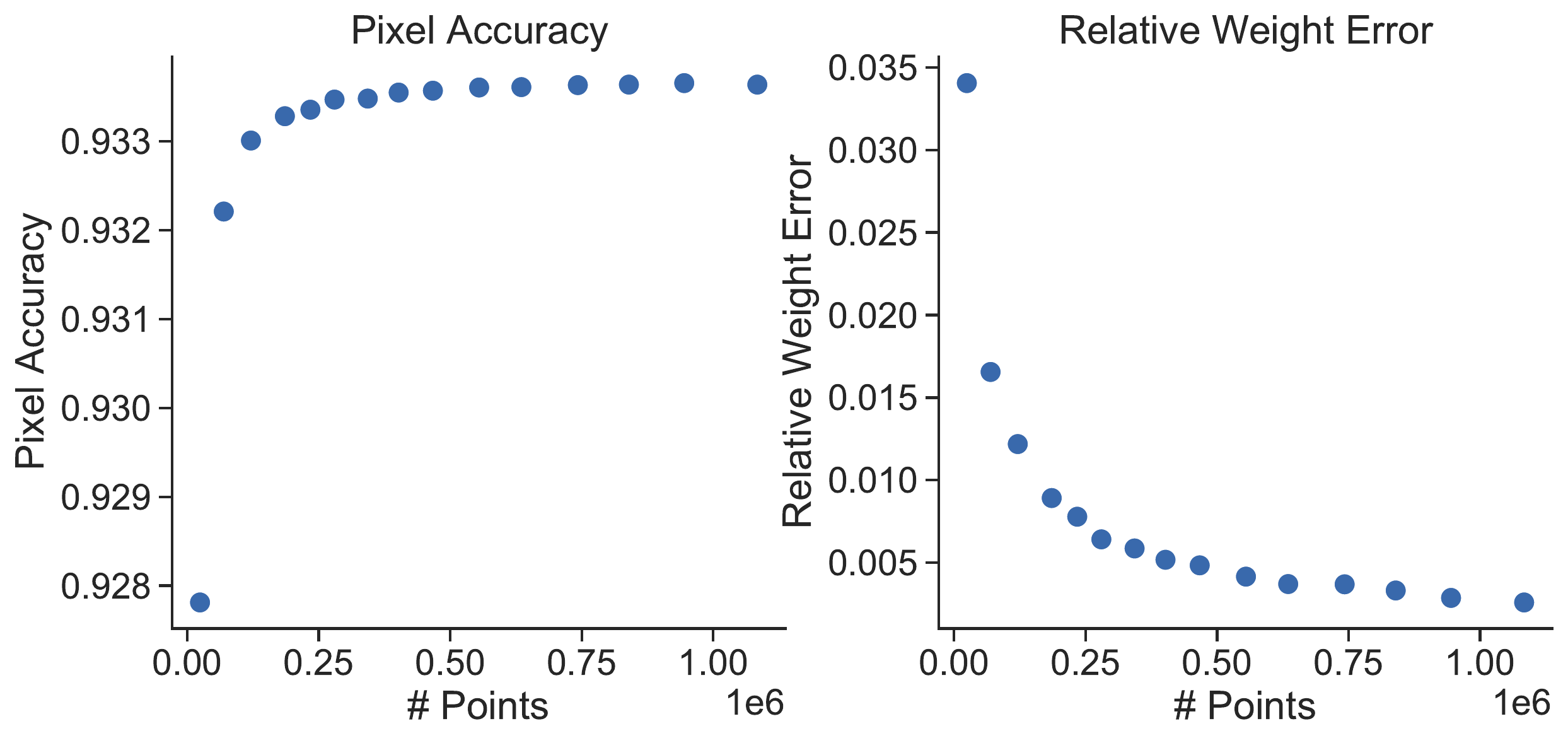}
	\caption{Accuracy (left) and relative cluster weight error (right) for different resolutions.}
	\label{fig:3dgrainscanresolutionwblsa}
\end{figure}

Again, we obtain an accuracy of $0.93$ for image supports above $100,000$ points. However, the resolution coreset seems to cause a relative cluster weight error which is only half of that for pencil coresets, i.e., $0.01$ instead of $0.025$ for the latter. With increasing resolution, $\Psi_G$ decreases further to $0.0022$ for 1 million points. \cref{fig:3dgrainscanresolutionwblsa} depicts the details.

The experimental results suggest that pertaining the grid structure in the image support yields, at the same level of sparsity, diagrams which are closer to the optimal diagram for the original instance. 

\subsubsection*{Additional grain scan-based interior removal}
\begin{figure}[h!]
	\centering
	\includegraphics[width=1\linewidth]{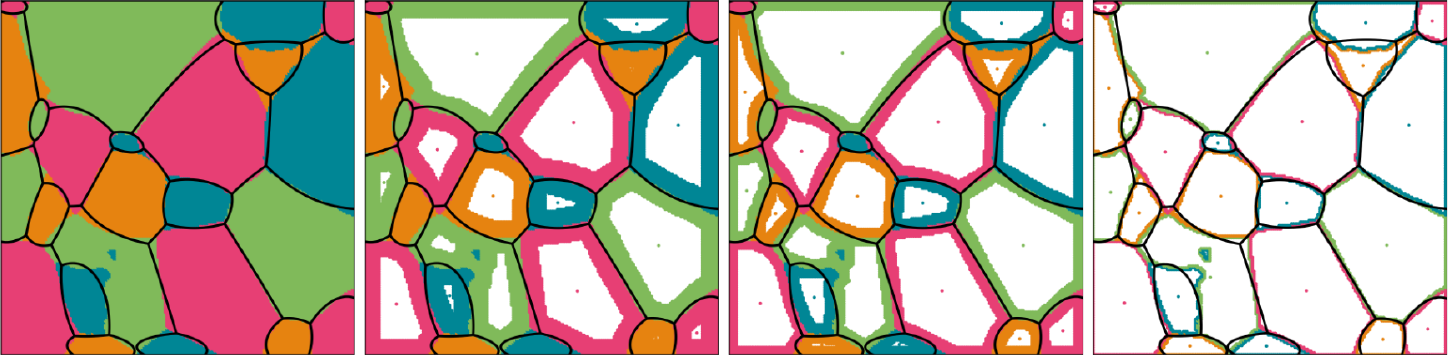}
	\caption{GBPDs from grain scan-based interior removal for decreasing $\bar\delta.$ From left to right: Computed diagrams (cell boundaries shown in black) based on decreasing image supports (the non-white and non-black colored pixels in the respective images).}
	\label{fig:grainscan-2d-interiorcoreset}
\end{figure}

In our practical evaluation of grain scan-based interior removal, we used a uniform setting, i.e., $\delta_i=\bar\delta$ for all $i\in [k]$ and employed the grid graph distance. See \cref{fig:grainscan-2d-interiorcoreset} for a visualization in 2D.

In order to be able to preserve the size constraints in \cref{lp:primal} for the reduced image support, we represent each removed $\bar\delta$-interior $G_i(\bar\delta)$ by a single point at $G_i$'s centroid, weighted by the cardinality of the removed interior points. While, particularly for large grains, a considerable number of points may be omitted, the total length of the unmodified \enquote{grain boundaries} is typically quite large. Hence, as a sole thinning strategy, interior removal is insufficient for reducing the image support to render the linear program tractable in practice. 

When combined with the previous coreset strategies, grain scan-based interior removal can, however, add a significant additional contribution. Pencil and resolution coresets still yield many points in each grain's interior, which can be removed for further speeding up $\cref{lp:primal}$.

\begin{figure}[h!]
	\centering
	\includegraphics[width=1\linewidth]{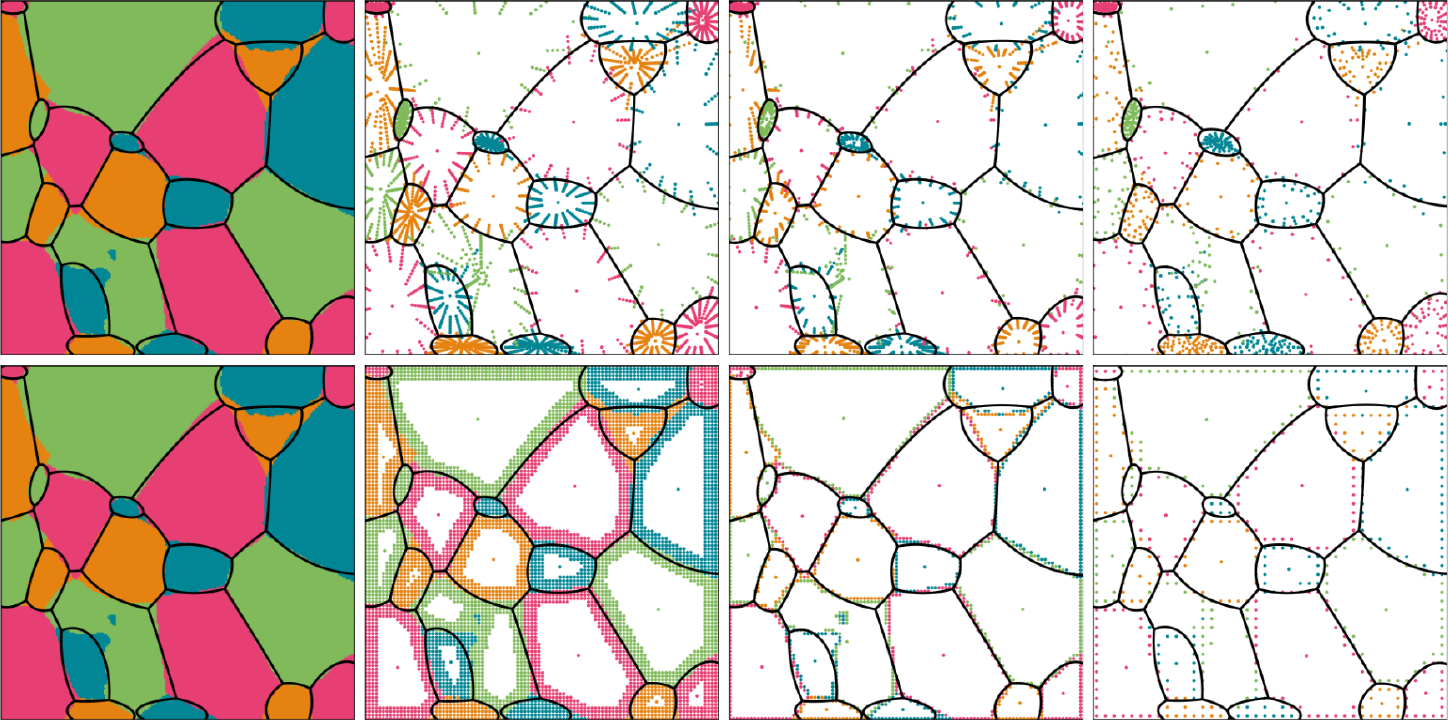}
	\caption{Pencil (top) and resolution coresets (bottom) combined with interior removal. From left to right: Computed diagrams (cell boundaries shown in black) based on decreasing image supports (the non-white and non-black colored pixels in the respective images).}
	\label{fig:grainscan-2d-interiorbatchingresolutioncoreset}
\end{figure}

 \Cref{fig:grainscan-2d-interiorbatchingresolutioncoreset} illustrates the combination of the coresets with interior removal in 2D. Again, the resulting diagrams are, as theory suggests, very similar to the original image.
 
In the experimental study, we investigated the effect of different parameter choices on the accuracy and the relative cluster weight error. Since we discussed the dependence of $\Phi_G$ and $\Psi_G$ on the coreset parameters already, we focus on the dependence on the parameter $\bar\delta$. Recall that the $\bar\delta$-interior $G_i(\bar\delta)$ decreases with increasing $\bar\delta$. Fortunately, already a grid graph distance of $2$ produced excellent fits. Best results were obtained for~$\bar\delta=4$. 

\begin{figure}[h!]
	\centering
	\includegraphics[width=1\linewidth]{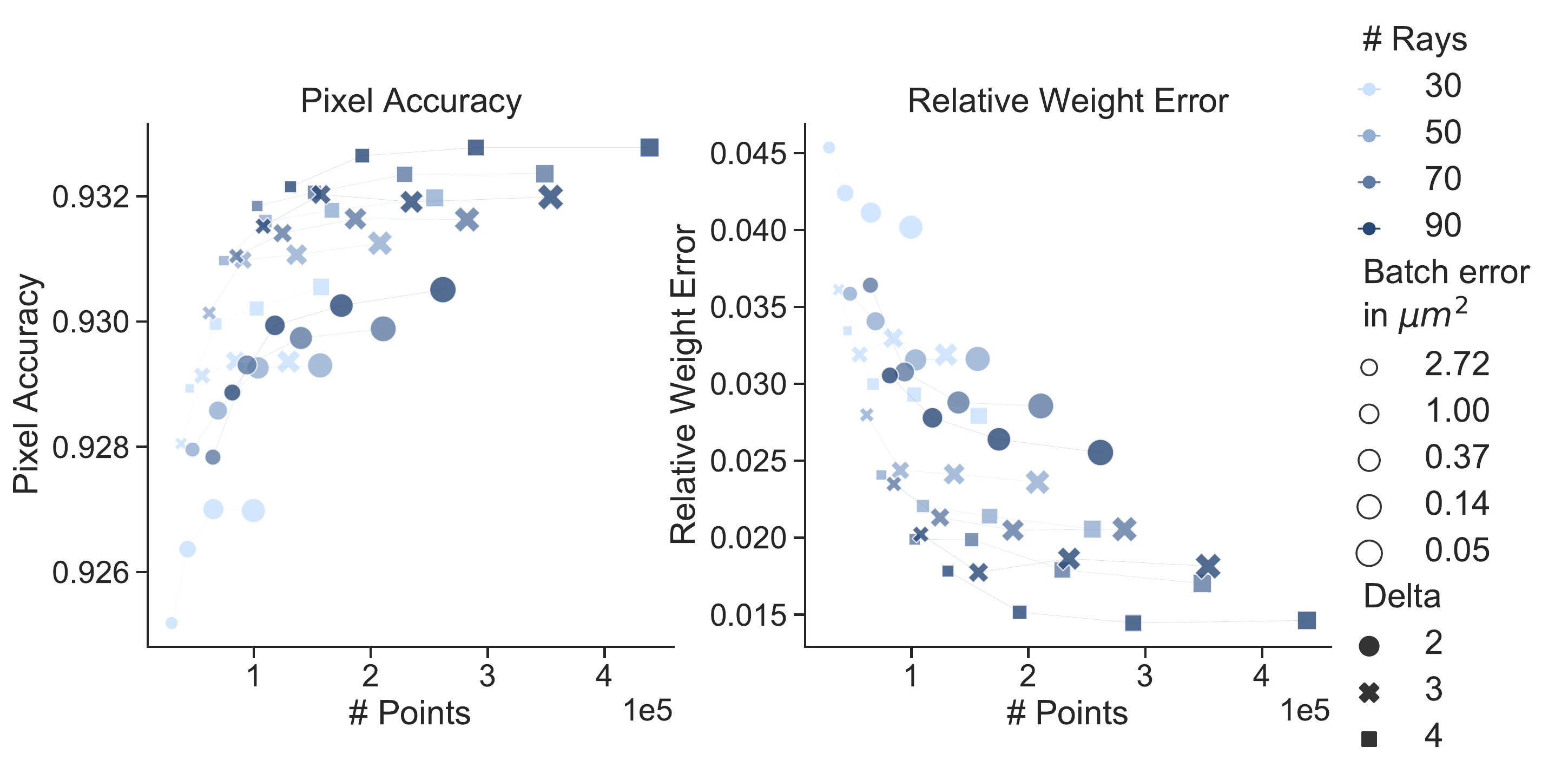}
	\caption{Combination of pencil coresets with interior removal: Dependence of $\Phi(G)$ and $\Psi_G$ on the size of the resulting image support. Again, the color indicates the number of rays, the size refers to the batch error while the shape of the symbol signifies the different values of $\bar \delta$.}
	\label{fig:3dgrainscanbatchinginteriorwblsa}
\end{figure}

\begin{figure}[h!]
	\centering
	\includegraphics[width=1\linewidth]{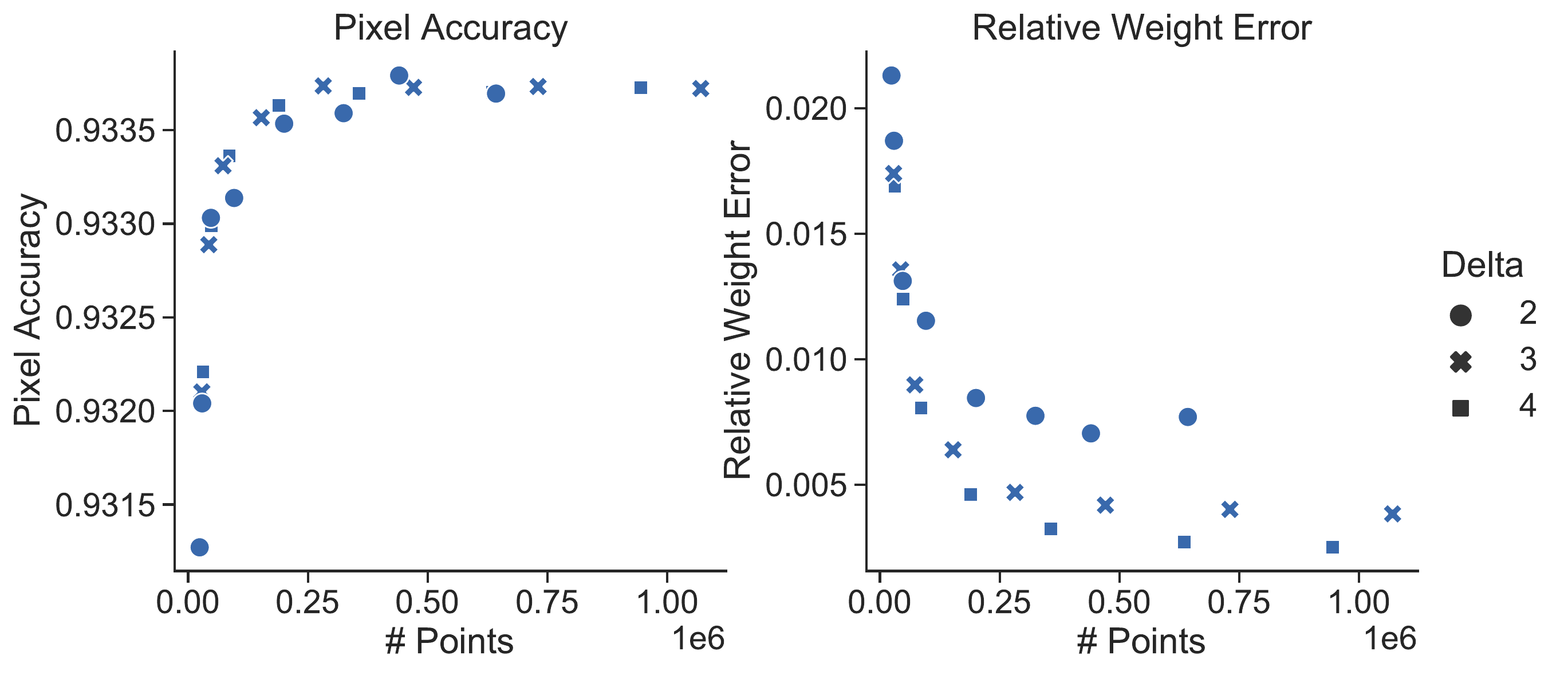}
	\caption{Combination of resolution coresets with interior removal: Dependence of $\Phi(G)$ and $\Psi_G$ on the size of the resulting image support.}
	\label{fig:3dgrainscanresolutioninteriorwblsa}
\end{figure}

Figures \ref{fig:3dgrainscanbatchinginteriorwblsa} and \ref{fig:3dgrainscanresolutioninteriorwblsa} show the effect of different choices for the parameters on the accuracy and the relative cluster weight error. 

In comparison, the combined thinning process based on resolution coresets outperformed that of pencil coresets: Its accuracy is higher, and its relative cluster weight error is smaller for all investigated coreset sizes. The difference is, however, generally only in the second or third digit. 

All in all, using as few as $85,828$ points was sufficient to obtain very high accuracy and low relative cluster weight error. Increasing the size of the image support to $1,000,000$ points results in only very marginal additional improvements. In our setting, an average of only $145$ points per grain sufficed for obtaining nearly the same quality of fit as the $116,476$ points in the original grain scan. 

\subsection{s-GBPD vs. \DiM} \label{sec:method_comparison}

While we showed in the previous paragraph that the indirect clustering-based methods yield already excellent fits for small sample sizes, we will now add a comparison with the direct benchmark method \DiM. Specifically, we computed $\apdmatrices, \apdsites, \apdweights$ via (\DiM) and derived the clustering $\hat{C}$ from the diagram $\apd(\apdmatrices, \apdsites, \apdweights)$. To decrease the number of constraints we reduced the the image support with the above techniques to $85,828$ points, and used only distance constraints for neighboring grains in the scan. 

In addition to the accuracy and the relative cluster weight error and, of course, the running time, we investigated various additional parameters. In particular, we computed the relative deviation of the centroids $c(G_i)$ and $C(\hat{C}_i)$, and also of the covariance matrices $A_i^{-1}= \Sigma_i=\text{Cov}(G_i)$ and $\text{Cov}(\hat{C}_i)$. More precisely, the \emph{relative centroid error} and the \emph{relative covariance error} are defined as
\begin{equation*}
\frac{1}{n} \sum_{i=1}^k \kappa_i\norm{c(G_i) - c(\hat{C}_i)}_{2} \quad \text{ and }
\quad \frac{1}{n}\sum_{i=1}^k \kappa_i\norm{\Sigma_i - \text{Cov}(\hat{C}_i)}_2.
\end{equation*}
In the latter case, the spectral norm measures the largest singular value of the covariance differences and, thus, essentially captures the difference in the largest singular vectors of the covariance matrices.

As characteristics of the grain topology affect the physical properties of polycrystals, we also address the cell combinatorics. In particular, we computed the percentage of \emph{correct neighborhoods} of the grains. More precisely, we count a voxel as adjacent to a given one, if it is contained in what is called the 26-neighborhood in discrete topology, i.e. if its centroid is among the 27 points which can be reached by at most one move in each coordinate direction of the grid graph. A grain is a neighbor of another grain if a voxel of the former is adjacent to a voxel of the latter.  The neighborhood is correct if all grain neighbors in the ground truth are also neighbors in the representation.

We also computed the percentage of grains whose neighborhood has at most to 1 or 2 errors, where an error is a missing or additional grain neighbor not present in the ground truth. In view of the results of the previous section, we focus on the resolution coresets with interior removal as the specific version of s-GBPD, with image support of 85,828 points. \Cref{tb:comparison_approximation_methods} depicts these results.

\begin{table}[h!]
	\centering
	\begin{tabular}{l|rr}
		\toprule
		\textbf{Performance characteristics}          	& \textbf{s-GBPD}   & \textbf{\DiM}    \\ \midrule
	  Accuracy $\Phi_G$ &  0.9334        & \textbf{0.9559} \\
Relative cluster weight error $\Psi_G$ & \textbf{0.01} & 0.02           \\
Relative centroid error (in $\mu$m) &  0.62         & \textbf{0.61}  \\
Relative covariance error (in $\mu$m${}^2$) & 15,53      & \textbf{8.72} \\
Correct neighborhoods (in \%) 	& 35.70 & \textbf{63.45} \\
Correct neighborhoods up to 1 error (in \%)  & 72.08 & \textbf{91.71} \\
Correct neighborhoods up to 2 errors (in \%) 	& 92.39 & \textbf{98.14} \\
Computation time (in min)       &  \textbf{2.60}        & 627.00          \\ \bottomrule
\end{tabular}
\caption{Comparison of s-GBPD and \DiM for polycrystal representation.}
\label{tb:comparison_approximation_methods}
\end{table}


While the Laguerre method from~\cite{Lyckegaard2011} achieves an accuracy~$\Phi_G$ of~0.86 on this data set (see~\cite{Lyckegaard2011}), both our methods relying only on parts of the data improve the accuracy to a value of over~0.93. As expected, the direct model algorithm \DiM performs better in most criteria but requires a significantly higher running time than s-GBPD. Since the differences are mostly marginal, the experimental findings indicate that we can compute a high precision 3D diagram representation of the grain scan in less than 3 minutes using s-GBPD. We expect that s-GBPD behaves similarly for other data sets as long as the resolution and the number of grains are of a similar magnitude.

\section{Outlook}\label{sec:conclusion}

As we have seen, s-GBPD generates high-quality diagram representations of polycrystalline scans at speeds acceptable in many applications. If still faster methods are required, heuristics will be the algorithms of choice; see \cite{AFGK22}. 

Let us close with some remarks on the role the grain scan plays in applying s-GBPD. While we focussed here on the situation that the whole grain scan $(X,Y)$ is available, s-GBPD can be adjusted to require only certain measurements of the polycrystal. Grain volumes, centers, and second-order moments can, in many cases, be measured directly using tomographic techniques, and the described clustering approach can be solely based on such measurements. Both, pencil and resolution coresets can still be constructed then. Grain scan-based interior removal, on the other hand, is (as emphasized correctly by its very name) based on the knowledge of $(X,Y)$. 

It is, however, possible to remove the dependency on $(X,Y)$ by adding an assumption that is common in  \emph{stochastic grain growth models}. In fact, in such models, the underlying grain growth is often viewed as some realization of a stochastic process involving multivariate normal distributions $\mathcal{N}(\mu_i,\Sigma_i)$ for each grain. Their $3$-dimensional mean vectors $\mu_i$ and $3\times 3$-covariance matrices $\Sigma_i$ are determined by the measurements and correspond to $\apdmatrices$ and $\apdsites$. Under such an assumption, interior points can be removed from the coreset by removing points with high density measured by the corresponding multivariate density function. As points close to the centroids have high density, this approach removes the interior points. Hence, it is possible to obtain additional speed-ups even when only the above measurements are available.

\section*{Acknowledgements}
We thank Henning Friis Poulsen for fruitful discussions and valuable comments on an earlier version of the present paper. We are also grateful to Henning Friis Poulsen and Allan Lyckegaard for providing the real-world data set used in this paper.

\printbibliography

\end{document}